\documentclass[12pt,twoside]{article}
\usepackage{amsmath,amsbsy,amsfonts}
\usepackage{hyperref}
\newtheorem{lem}{Lemma}[section]
\newtheorem{prop}{Proposition}[section]
\newtheorem{thm}{Theorem}[section]

\newtheorem{rmk}{Remark}[section]
\newtheorem{claim}{Claim}[section]
\newtheorem{defi}{Definition}[section]

\let\Section=\section

\def\section{\setcounter{equation}{0}\Section}

\def\nd{\noindent}

\def\hsf{\hspace*{\fill}}

\def\proof{{\rm \bf Proof}} 

\newcommand{\w}{W_0^{1,\Phi}(\Omega)}

\def\r{\bf{R}}

\begin{document}
\title{ Multiple Positive Solutions for a Class of   \\
Nonlinear Elliptic Eigenvalue Problems\\
with  a Sign-Changing Nonlinearity\thanks{Partially supported by PROCAD/CAPES/UFG/UnB-Brazil}}

\author{M. L. Carvalho\thanks{M. L. Carvalho was partially supported by CAPES/Brazil} ~~   J. V. Goncalves\thanks{J. V. Goncalves was partially supported by CNPq/Brazil}\,\, ~  K. O. Silva\thanks{K. O. Silva was partially supported by CAPES/Brazil}  }

\date{}

\pretolerance10000

\maketitle

\begin{abstract}
\noindent {\small In 2009 Loc and Schmitt   established a result on sufficient conditions for multiplicity of solutions of a class of nonlinear eignvalue problems for  the p-Laplace operator under  Dirichlet boundary conditions, extending an earlier result of 1981  by Peter Hess  for the Laplacian. Results on necessary conditions for existence were also established. In the present paper the authors extend the main results by Loc and Schmitt   to  the $\Phi$-Laplacian. To overcome the difficulties with this much more general operator it was necessary to employ regularity results by Lieberman, a strong maximum principle by Pucci and Serrin and a general result on lower and upper solutions by Le \cite{Khoi}.}
\vskip.5cm

{\small 

}
\end{abstract}

\section{Introduction}

\nd  We study the nonlinear eigenvlaue problem
~\begin{equation}\label{3.1}
\left\{\begin{array}{rllr}
-{\rm div} \big( \phi(|\nabla u|)\nabla u\ \big ) &=& \lambda f(u)~~\mbox{in}~~\Omega,\\

u&=&0~~ \mbox{on}~~  \partial {\Omega,}
 \end{array}
\right.
\end{equation}
\nd where $\Omega\subset\mathbb{R}^N$ is  bounded domain with smooth boundary $\partial \Omega$,~$\lambda > 0$ is a parameter and $\phi : (0, \infty) \to (0, \infty)$ is a  $C^1$-function satisfying 
\begin{itemize}
  \item[($\phi_1$)]  \     \ $\mbox{(i)} \ \ t\phi(t)\to 0 \ \mbox{as} \  t\to 0, \\  
\\
 ~~ \mbox{(ii)} \  t\phi(t)\to\infty \ \mbox{as} \  t\to\infty$, 
  \item[($\phi_2$)]  \ $t\phi(t) \ \mbox{\it is strictly increasing in}~ (0, \infty)$,
 \end{itemize}
\nd while $f : [0, \infty) \to [0, \infty)$ is a continuous function satisfying

\begin{itemize}
  \item[($f_1$)] \  $f(0)\geq 0$,
  \item[($f_2$)] \  {\it there exist positive  numbers} $a_k, b_k,\ k = 1, \cdots, m$ {\it such that}
$$
0<a_1<b_1<a_2<b_2<...<b_{m-1}<a_m, \\
$$
$$
f(s)\leq 0 \ \mbox{ if } \ s \in (a_k,b_k), \\
$$
$$
f(s)\geq 0 \ \mbox{ if } \ s \in (b_k,a_{k+1}),
$$
\item[($f_3$)] \  $\displaystyle\int_{a_{k}}^{a_{k+1}}f(s)ds>0$,  $k = 1,  \cdots,  m-1$.
\nd
\end{itemize}

\begin{rmk}
We  extend   $t \mapsto t \phi(t)$ to the whole of $\r$ as an odd function.
\end{rmk}

\nd In \cite{Hess}, by means of vatiational and topological methods and arguments with lower and upper solutions, Hess proved a result on existence of multiple solutions of   (1.1) for the case of the  Laplacian operator which means taking $\phi(t)  \equiv 1$  in problem (1.1). 
\vskip.1cm

\nd  According to Hess, the results in \cite{Hess} were motivated by Brown \& Budin \cite{BB-1, BB-2}
which in turn were motivated by  the literature on nonlinear heat generation.
\vskip.1cm

\nd In  \cite{loc-schmitt}, Loc \& Schmitt extended the result by Hess to the $p$-Laplacian operator by means  of taking $\phi(t) = t^{p-2} \ \mbox{with} \  1 < p < \infty$ in (1.1). Actually, in \cite{loc-schmitt}  the authors showed (see also  Dancer \& Schmitt \cite{dancer-schmitt}) that $(f_1)-(f_3)$ are suffitient conditions for the existence of $m-1$ solutions of
$$
- \Delta_p  u = \lambda f(u)~\mbox{in}~\Omega,~~u = 0~\mbox{on}~\partial \Omega, 
$$
\nd for $\lambda$ large,  while   if $(f_1)-(f_2)$ hold, and $u$ is a solution of the problem above 
with $a_k < \|u\|_\infty \leq a_{k+1}$ then $(f_3)$ holds.
\vskip.3cm

\nd In the present  work  we were able to adapt the techniques in  Hess \cite{Hess} and in Loc \& Schmitt \cite{loc-schmitt}  to prove a result extending the main theorem of \cite{loc-schmitt}  to the more general operator 
$$
 \Delta_{\Phi}  u  :={\rm div} \big( \phi(|\nabla u|)\nabla u\ \big ) 
$$
\nd named $\Phi$-Laplacian, where $\Phi$ is the even function  defined by
$$
\Phi(t) = \int_0^t s \phi(s) ds,~~ t \in \r.
$$
\nd and the function $\phi$ satisfies $(\phi_1)-(\phi_2)$ and the following further conditions
$$
\begin{array}{lrl}
(\phi_3)~~~~   \displaystyle \sum_{i,j=1}^N\frac{\partial \alpha_j}{\partial\eta_i}(\eta) \xi_i \xi_j\geq\Gamma_1~ \phi(|\eta|)|\xi|^2,\\

(\phi_4)~~~ ~~  \displaystyle |\sum_{i,j=1}^N\frac{\partial \alpha_j}{\partial\eta_i}(\eta)|\leq\Gamma_2~ \phi(|\eta|),
\end{array}
$$
\nd where  $\Gamma_1,\Gamma_2 > 0$ are  constants,
$$
\xi = (\xi_1, \cdots, \xi_N),~~\eta=(\eta_1,...,\eta_N),
$$
$$
\alpha_j(\eta)=\phi(|\eta|)\eta_j,~~j=1, \cdots, N.
$$ 
\nd It is well known that the $p$-Laplacian is included in this class of operators. Moreover, our result includes a broader class of operators, for example $\Delta_{\Phi}$ with 
\begin{equation}\label{ex 1}
\Phi(t)=(1+t^2)^{\gamma} - 1~~ \mbox{where}~~\gamma>\frac{1}{2}
\end{equation} 
\nd and
\begin{equation}\label {ex 2}
\Phi(t)=t^p \log(1+t)~~ \mbox{where}~~p\geq 1.
\end{equation}
\nd See the Appendix for further comments on these examples. 

 \begin{defi}\label{1.0}  
By a  solution of   $(\ref{3.1})$  we mean a function $u \in C_0^{1}(\overline{\Omega})$  satisfying
$$
\int_\Omega\phi(|\nabla u|)\nabla u \cdot \nabla v dx =\lambda\int_\Omega f(u)v dx,\ \  v\in C_0^{1}(\overline{\Omega}),
$$
\nd where 
$$
C_0^{1}(\overline{\Omega}) = \{ u \in C^{1}(\overline{\Omega})~|~u = 0~ \mbox{on}~ \partial \Omega \}.
$$
\end{defi}

\nd Our main result below extends Theorem 1.1 by Loc \& Schmitt in \cite{loc-schmitt} to the more general operator $\Delta_{\Phi}$.

\begin{thm}\label{1.1} 
Assume  $(\phi_1)-(\phi_4)$.  Then 
\vskip.1cm

\nd {\rm (i)}~ if $(f_1)-(f_3)$ hold, there is $\overline{\lambda}>0$, such that for each $\lambda>\overline{\lambda}$,  $(\ref{1.1})$ admits at least $m-1$ solutions say $u_1,...,u_{m-1}$  such that 
$$
a_1 <  \|u_{1}\|_\infty \leq a_2 < \|u_{2}\|_\infty \leq \cdots \leq a_{m-1}<\|u_{m-1}\|_\infty \leq a_{m},
$$
\nd {\rm (ii)}~  if $u$ is a solution of $(\ref{3.1})$ with $a_k < \|u\|_\infty \leq a_{k+1}$ and $(f_1)-(f_2)$  hold then  $(f_3)$  also  holds.
\end{thm}

\nd Due to the more general nature of  $\Delta_{\Phi}$, in our proof of  theorem  \ref{1.1}  above  it was necessary to get into the framework of Orlicz-Sobolev spaces. It was also necessary to employ  regularity results by Lieberman \cite{Lieberman, lieberman2}, a  strong maximum principle due to Pucci \& Serrin \cite{Pucci} which holds in our setting as well as a more general result on lower and upper solutions due to  Le \cite{Khoi}.  

\section{Notations and Auxiliary Results}

\nd Consider the family of problems associated to (\ref{3.1})
~\begin{equation}\label{3.1k}
\left\{\begin{array}{rllr}
-\Delta_{\Phi} u &=& \lambda f_k(u)~~\mbox{in}~~\Omega,\\

u&=&0~~ \mbox{on}~~  \partial {\Omega,}
 \end{array}
\right.
\end{equation}
\nd where for   $k = 2, \cdots, m.$, $f_k : \r \to \r$ is the continuous function 
$$
f_k(s) = \left\{ \begin{array}{rl}
 f(0) &\mbox{if $s\leq0$} ,\\
  f(s) &\mbox{if $0\leq s\leq a_k$ }, \\
  0 &\mbox{if $s>a_k$}.
       \end{array} \right.
$$

\nd  In this work  
$$
 W^{1, \Phi}(\Omega) =\Big\{u \in L_\Phi(\Omega)~|~ \frac{\partial u}{\partial x_i} \in L_\Phi(\Omega),~ i=1,...,N \Big\},
$$
\nd is the Orlicz-Sobolev space, where $L_\Phi(\Omega)$ is the Orlicz space defined through the $N$-function $\Phi$, endowed with the (Luxembourg) norm
$$
\|u\|_\Phi=\inf\left\{\lambda>0~|~\int_\Omega \Phi\left(\frac{u(x)}{\lambda}\right) dx \leq 1\right\},
$$
\nd while  $\w$ denotes  the closure of  ${C}^\infty_0(\Omega)$ with respect to the usual norm of $W^{1, \Phi}(\Omega)$. We refer the reader to Adams \cite{A}, concerning Orlicz-Sobolev spaces.
\begin{rmk} \label{obss} {\rm The reader is referred to $\cite{A}$, $\cite{Fukagai}$ for the basic results below:
\vskip.1cm

\nd {\bf (i)}  if $\phi$ satisfies~  $(\phi_1)-(\phi_2)$~ it is an easy matter to check that  $\Phi$ is an $N$-function (or Young function),
\vskip.2cm

\nd {\bf (ii)}  if $\phi$ satisfies~ $(\phi_3)-(\phi_4)$~then, (cf. proposition  \ref{equivalencia-phi34} in the Appendix),  
   \begin{equation}\label{relacao a_1- a_2}
        \Gamma_{1}\leq \frac{(t\phi(t))'}{\phi(t)}\leq \Gamma_{2},~~t > 0,
   \end{equation} 
\vskip.2cm

\nd {\bf (iii)}   if  (\ref{relacao a_1- a_2}) holds then (cf. remark \ref{Aa} in the Appendix) there exist constants $\gamma_1,\gamma_2>1$  such that 
\begin{equation}\label{Delta_2}
  \gamma_1\leq\frac{t\Phi'(t)}{\Phi(t)}\leq \gamma_2,\  t>0,
\end{equation}
\vskip.2cm

\nd {\bf (iv)} It follows by  \cite[pg 542]{Fukagai}  and \cite[thm 8.20 pg 274]{A} that 
$L_\Phi(\Omega)$ is reflexive if condition (\ref{Delta_2})  holds true.
\vskip.2cm

\nd {\bf (v)} As a consequence of the remarks  {\bf (i)-(iv)}  above, the spaces $L_\Phi(\Omega)$ and 
$W^{1, \Phi}(\Omega)$ are reflexive if $\phi$ satisfies $(\phi_1)-(\phi_4).$ }
\end{rmk}
\vskip.2cm

\nd The energy functional  associated to $(\ref{3.1k})$ is
$$
I_k(\lambda,u)=\int_\Omega\Phi(|\nabla u|) dx-\lambda\int_\Omega F_k(u) dx,~ u \in \w,
$$
\nd where 
$$
F_k(s)=\int_0^s f_k(t)dt.
$$
\nd  It is known  that  
$I_k(\lambda,\cdot):W_0^{1,\Phi}(\Omega)\to\mathbb{R}$ is a $C^1$-functional and
$$
\displaystyle \langle  I_k^{\prime}(\lambda,u), v \rangle =  \int_\Omega\phi(|\nabla u|)\nabla u \cdot \nabla v   dx - \lambda\int_\Omega f_k (u)v dx,\  v\in W_0^{1, \Phi}(\Omega).
$$
\nd Thus,  a critical point $u \in \w$  of $I_k(\lambda,\cdot)$ is a weak solution of $(\ref{3.1k})$,  in the sense that 
$$
\int_\Omega\phi(|\nabla u|)\nabla u \cdot \nabla v   dx = \lambda\int_\Omega f_k (u)v dx,\  v\in W_0^{1, \Phi}(\Omega).
$$
\begin{rmk}\label{regularity}
If $u$ is a weak solution of  $(\ref{3.1k})$  then, since $f_k$ is bounded and continuous, $f_k(u) \in L^{\infty}(\Omega)$. It follows by Lieberman $\cite[theorem~ 1.7]{Lieberman}$  that $u \in C^{1,\alpha}({\overline{\Omega}})$ where $\alpha \in (0,1)$ and so $u$ is a solution of  $(\ref{3.1k}) $ in the sense of definition  $\ref{1.0}$. 
\end{rmk}

\section{Technical Lemmata}

\nd The result below is crucial in this paper, it was proved by Loc \& Schmitt for Sobolev spaces and its proof in our case is similar. We leave its  proof to the end of the section.
\vskip.2cm

\begin{lem}\label{lema 2.1}
 Let $g:\mathbb{R}\to\mathbb{R}$ be a continuous function such  that $g(s)\geq 0$ for $s\in (-\infty,0)$ and assume that there is some $s_0\geq 0$  such that $g(s)\leq 0$ for $s\geq s_0$.  Let $u \in \w$ be a weak solution of

\begin{equation}\label{2.1}
 \left\{ \begin{array}{c}
 -\Delta_\Phi u= g(u) ~\mbox{ in}~~ \Omega, \\
u   = 0~\mbox{on}~\partial \Omega .
       \end{array} \right.
\end{equation}
\nd Then   $0 \leq u \leq s_0~\mbox{a.e. in}~\Omega$.
\end{lem}
\begin{lem}\label{Qk not-empty}
Let $\lambda > 0$. Then there is   $v_k \equiv v_k(\lambda) \in \w$ such that
$$
I_{k}(\lambda, v_k) = \min_{u \in \w } I_{k}(\lambda, u).
$$
\end{lem}

\nd   \proof . It is enough to show that $I_{k}(\lambda, \cdot)$ is both coercive and weakly sequentially lower semicontinuous, (w.s.l.s.c.  for short).
\vskip.2cm

\nd To show the coerciveness, by  lemma \ref{lema_naru}, (cf. Appendix),  the continuous embedding $\w \hookrightarrow L^{1}(\Omega)$ and the Poincar\'e inequality, (cf \cite{Gz1}),  we have
\begin{eqnarray}
  I_k(\lambda,u) &\geq& \min\{ \|\nabla u\|^{\gamma_1},\|\nabla u\|^{\gamma_2} \}-\lambda C |u|_1     \nonumber \\
   &\geq& \min\{\|\nabla u\|^{\gamma_1},\|\nabla u\|^{\gamma_2}\}-\lambda C\|\nabla u\|, \nonumber
\end{eqnarray}
\nd which shows that $I_{k}(\lambda, \cdot) : \w \to \r$  is coercive.
\vskip.2cm

\nd To show that $I_{k}(\lambda, \cdot)$ is w.s.l.s.c, at first notice that 
$$
u \in \w \mapsto \int_\Omega \Phi(|\nabla u|) dx
$$
\nd is continuous and convex. Take $(u_n)$ such that  $u_n\rightharpoonup u$ in $W_0^{1,\Phi}(\Omega)$.  Using the embedding  (cf. Adams \cite{A}),
$$
\w \stackrel {\scriptsize cpt}   \hookrightarrow  L_{\Phi}(\Omega) ,
$$
\nd  and arguments with the convexity of $\Phi$, there is $h \in L_{\Phi}(\Omega)$ such that
$$
u_n \rightharpoonup u~\mbox{in}~~ \w,~~ u_n \to u~ \mbox{and}~~  |u_n|\leq h~~ \mbox{a.e. in}~\Omega.
$$
\nd By Lebesgue's Theorem,
$$
\int_\Omega F_k(u_n) dx \to \int_\Omega F_k(u) dx.
$$
\nd   Since $\Phi$ is continuous and convex,
$$
\int_\Omega \Phi(|\nabla u|) dx \leq \liminf \int_\Omega \Phi(|\nabla u_n|) dx.
$$
\nd It follows that 
$$
I_{k}(\lambda, u) \leq \liminf I_{k}(\lambda, u_n).
$$
\nd As a consequence, there is minimum $v_k \equiv v_k(\lambda)$ of $I_{k}(\lambda, \cdot)$. \hfill \fbox \hsf
\vskip.2cm

\begin{lem}\label {lema 2.3}
  There is $\lambda_k>0$ such that  
$$
a_{k-1}<\|v_{k}\|_\infty \leq a_{k}
$$
\nd for each minimum $v_k \equiv v_k(\lambda)$  of $I_{k}(\lambda, \cdot) $ with $\lambda > \lambda_k$.
\end{lem}

\nd \proof\ {\bf of Lemma \ref{lema 2.3}} The proof is similar to  the ones   in \cite{Hess, loc-schmitt}. So   we will just sketch the main steps. Take $\delta > 0$ and consider the open set
 $$
\Omega_\delta=\{x\in\Omega~|~ \operatorname{dist}(x,\partial\Omega)<\delta\}.
$$ 
\nd Set  
$$
{\widetilde{\alpha }}_k :=F(a_k)-\max{\{F(s)~|~ 0\leq s\leq a_{k-1}\}}
$$
\nd and note that by $(f_3)$, $ {\widetilde{\alpha }}_k > 0$. Choose  $w_\delta\in C_0^\infty(\Omega)$ such that 
$$
0\leq w_\delta\leq a_k~\mbox{and}~ w_\delta(x)=a_k,~  x\in\Omega\setminus\Omega_\delta.
$$
\nd Writing $\Omega=\Omega_{\delta} \cup (\Omega \setminus \Omega_{\delta})$ and setting $ C_k=\max{\{|F(s)|~|~\ 0\leq s\leq a_k\}}$ we get to,

$$
 \int_\Omega F(w_\delta) dx \geq \int_\Omega F(a_k) dx -2C_k|\Omega_\delta|. \nonumber
$$

\nd Let  $u\in W_0^{1,\Phi}(\Omega)$ such that $0\leq u\leq a_{k-1}$. By the  inequality above we have
$$
\int_\Omega F(w_\delta) dx -\int_\Omega F(u) dx \geq {\widetilde{\alpha }}_k |\Omega|-2C_k|\Omega_\delta|.
$$
\nd Since  $|\Omega_\delta|\to 0$ as $\delta\to 0$ there is $\delta>0$ such that 
$$  \eta_k := {\widetilde{\alpha }}_k |\Omega|-2C_k|\Omega_\delta|>0.
$$
\nd Set $w=w_\delta$ and pick $u\in W_0^{1,\Phi}(\Omega)$ with $0\leq u\leq a_{k-1}$. Choosing $\lambda_k > 0$ large enough, taking $\lambda \geq \lambda_k$ and  making use of the expessions of  $I_k(\lambda,w), I_{k-1}(\lambda,u)$ and the inequality just above we infer that
\begin{equation}
 I_k(\lambda,w)-I_{k-1}(\lambda,u)  \leq \int_\Omega\Phi(|\nabla w|) dx -\lambda\eta_k 
   < 0
\end{equation}
\nd and  hence  
\begin{equation}\label{1.4}I
_k(\lambda,w)<I_{k-1}(\lambda,u)~\mbox{ for }~\lambda \geq \lambda_k. 
\end{equation}

\nd To finish, assume, on  the contrary,  that there is a minimum $v_k(\lambda)$ of  $I_k(\lambda,\cdot)$  such that $v_k(\lambda) \leq a_{k-1}$. It follows by $(\ref{1.4})$ and lemma \ref{lema 2.1} that 
$$
I_k(\lambda,w)<I_{k-1}(\lambda,v_k(\lambda)).
$$
\nd   On the other hand, since 
$v_k(\lambda)$ is a minimum of $I_k(\lambda,\cdot)$ we have
$$
 I_k(\lambda,v_k)  \leq  I_k(\lambda,w)      \nonumber \\
  $$
\nd  The definitions of $I_k(\lambda, \cdot)$ and $I_{k-1}(\lambda, \cdot)$ and the inequalities just above  lead   to a contradiction.
\nd  This ends the proof of lemma \ref{lema 2.3}. \hfill \fbox \hsf
\vskip.3cm

\nd {\bf Proof of Lemma \ref{lema 2.1}}~ Let  $u \in W_0^{1,\Phi}(\Omega)$ be a weak solution of (\ref{2.1}).  Recall   that (even for Orlicz-Sobolev spaces)    $u^-=\max\{-u,0\} \in W_0^{1,\Phi}(\Omega)$. We have
\begin{equation}\label{No. 1}
 \int_\Omega \phi(|\nabla u|)\nabla u\nabla u^- dx  = \int_\Omega g(u)u^-  dx      
   = -\int_{u<0}g(u)u dx \geq 0. 
\end{equation}
\nd Moreover (also for  Orlicz-Sobolev spaces) one has
$$
\nabla u^- =
 -\nabla u~  \chi_{\{u < 0 \}}~\mbox{a.e. in }~\Omega.
$$
\nd Using this in   (\ref{No. 1})    we find that
 $$
\int_{\Omega} \phi(|\nabla u|)|\nabla u|^2  \chi_{\{u < 0 \}} dx =0
$$
\nd which shows that   $u\geq 0$ in $\Omega$.
\vskip.2cm

\nd  Let $(u-s_0)^+=\max\{u-s_0,0\}$. By using an argument as the one  above, we infer that $u\leq s_0$. This proves lemma \ref{lema 2.1}. \hfill \fbox \hsf

\section{Proof of Theorem \ref{1.1} } The proof is based on Loc \& Schmitt \cite{loc-schmitt}.  However, we will get into  datails   taking into account the Orlicz-Sobolev  spaces framework. In this sense we will make use of a general result on lower and upper solutions by Le \cite[theorem 3.2]{Khoi}  and a general strong maximum principle  by Pucci \& Serrin \cite[theorem 1.1]{Pucci}  which hold in our setting.
\vskip.2cm

\nd The proof of {\rm(i)} is easier. For the proof of {\rm(ii)} we will need the two lemmas below whose proofs  are left to the end of this section.
\vskip.2cm

\begin{lem}\label{lemma 3.1 Kaye}
Let $u$ be a  solution of $(\ref{1.1})$ in the sense of definition $\ref{1.0}$.
If  $u \geq 0$ and $f(0)>0$ then $u>0$ in $\Omega$.
\end{lem}

\nd The proof of this lemma, which is left to the end of this section, strongly uses a general form of the maximum principle due to Pucci \&  Serrin \cite{Pucci}.
\vskip.3cm

\nd In order to state the second lemma take an open ball $B$ centered at $0$ with radius $R$ containing $\Omega$. Consider the functions  $\alpha, \beta : {\overline{B}} \to \r$ defined as follows:
$$
 \alpha(x)= \left\{ \begin{array}{rll}
 u(x), &~ x \in  \overline{\Omega}  \\
  0, &~  x \in \overline{B} \setminus \Omega,
       \end{array} \right.
~~~~~ \beta(x)=a_{k+1},~   x \in \overline{B}.
$$
\nd Since $u \in  W_0^{1,\Phi}(\Omega) $ there is a sequence $\{u_n \} \subseteq  C_0^{\infty}(\Omega)$ such that  $u_n \rightarrow u$ in the norm of $W_0^{1,\Phi}(\Omega)$. Extending each $u_n$ to $B$ as zero outside $supp(u_n) \subset \Omega$ it follows that   $\{u_n \} \subseteq  C_0^{\infty}(B)$ and  $u_n \rightarrow \alpha$ in the norm of $W_0^{1,\Phi}(B)$ so that $\alpha \in  W_0^{1,\Phi}(B) $. 

\begin{lem}\label{lemma 3.2 Kaye}
$\beta$ and $\alpha$ are respectively upper and lower solutions of 
\begin{equation}\label{3.2}
 \left\{ \begin{array}{cl}
 -\Delta_\Phi u=\lambda f(u)~\mbox{in}~ B, \\
  \\
    u\in W_0^{1,\Phi}(B).
       \end{array} \right.
\end{equation}
\end{lem}

\nd The proof of this lemma, which is left to the end of this section,  uses a general theorem on lower and upper solutions due to Le \cite{Khoi}.
\vskip.3cm

\nd {\bf Proof of  {\rm(i)} of theorem \ref{1.1}}. Take  $\lambda > 0$. By lemma  \ref{Qk not-empty},   for each  $k = 2, \cdots, m$ there is a minimum $v_k \equiv v_k(\lambda)$ of $I_k(\lambda)$,  which is actually a weak solution of problem (\ref{3.1k}).  By remark  \ref{regularity},   $v_k \in C^{1}(\overline{\Omega})$  and by lemma \ref{lema 2.1}, $0 \leq v_k \leq a_{k}~\mbox{a.e. in}~\Omega$.
\vskip.2cm

\nd By lemma \ref {lema 2.3}, there is $\overline{\lambda} \geq \displaystyle \max_{2 \leq k \leq m} \{\lambda_k \}$ such that for $\lambda > \overline{\lambda}$, $v_2, \cdots, v_m$ are  solutions of problem (\ref{3.1}) satisfying
$$
a_1 <  \|v_{2}\|_\infty \leq a_2 < \|v_{3}\|_\infty \leq \cdots \leq a_{m-1}<\|v_{m}\|_\infty \leq a_{m}
$$
\nd We set $u_{k-1}  \equiv v_{k}(\lambda),~k = 2, \cdots, m$. This ends the proof of the first part of theorem \ref{1.1} .
\vskip.2cm

\nd {\bf Proof of  {\rm(ii)} of theorem \ref{1.1}}. We distinguish between two cases.
\vskip.3cm

\nd {\bf Case 1}~~ $f(0) > 0$.
\vskip.2cm

\nd This case is  more difficult. In order to address it we state and prove the lemma below.
\begin{lem}\label{thm 3.1 Kaye}
Assume $(\phi_1)-(\phi_4)$, $(f_1)-(f_2)$ and $f(0)>0$. If $u$ is a non-negative weak solution of $(\ref{1.1})$ such that  $a_{k-1} < \|u\|_\infty \leq a_{k}$ then
 $$\int_{a_k}^{a_{k+1}}f(s)>0.$$
\end{lem}

\nd \proof~  Let us think of $k = 2$, for a while. Take the lower and upper solutions respectively $\alpha$ and $a_2$ of  (\ref{3.2}). 
\vskip.1cm

\nd Applying theorem 3.2 of \cite{Khoi}  there is a  maximal solution say $\overline{u}$ of  (\ref{3.2}) such that $\alpha (x)\leq\overline{u}(x)\leq a_{2}$ for $x\in B$.
\vskip.2cm

\nd By remark  \ref{regularity},  $\overline{u} \in C_{0}^{1}(\overline{B})$ and by lemma \ref{lemma 3.1 Kaye}, $\overline{u} > 0~\mbox{in}~B$.
 \begin{claim}\label{rad symm}
 {\it $\overline{u}$ is radially symmetric, e.g.  $\overline{u}(x_1)=\overline{u}(x_2),~x_i \in B,~ |x_1|=|x_2|$}.
\end{claim}
\nd Indeed, assume on the contrary that 
$$
\overline{u}(x_1)<\overline{u}(x_2)~ \mbox{ for some}~ x_1,x_2 \in B~\mbox{ with}~ |x_1|=|x_2|.
$$
\vskip.2cm

\nd Choose a rotation matrix $P$ such that $x_2=P x_1$. Recall that $P^\top P=I$ and $|\operatorname{det}P|=1$. 
\vskip.1cm

\nd Set $u_1(x)=\overline{u}(Px)$. Since 
 $$
\nabla u_1(x)=P\nabla \overline{u}(Px),~  x\in\Omega, 
$$
\nd it follows that,  ($P$ is an isometry), 
$$
|\nabla u_1(x)| = |\nabla\overline{u}(Px)|. 
$$
\nd We contend that  $u_1$ is a weak solution of $(\ref{3.2})$. Indeed,  let $\varphi\in W_0^{1,\Phi}(B)$ and set $\psi(x)=\varphi(P^\top x)\in W_0^{1,\Phi}(B)$. We have by easy computation,

\begin{eqnarray}
  \int_B \phi(|\nabla u_1(x)|)\nabla u_1(x)\nabla \varphi(x) dx &=& \int_B \phi(|\nabla \overline{u}(Px)|)\nabla \overline{u}(Px)\nabla \varphi(x)   dx    \nonumber \\
   &=& \int_B \phi(|\nabla \overline{u}(y)|)\nabla \overline{u}(y)\nabla \psi(y)|\operatorname{det}(P)| dy
    \nonumber \\
   &=& \lambda\int_B f(\overline{u}(y))\psi(y) dy\nonumber \\
   &=& \lambda\int_B f(\overline{u}(Px))\psi(Px)|\operatorname{det}P^\top| dx\nonumber \\
   &=& \lambda\int_B f(u_1(x))\varphi(x) dx,\nonumber
\end{eqnarray}
\nd showing that $u_1$ is a solution of $(\ref{3.2})$. 
\vskip.2cm

\nd Of course,  $u_1$ is a subsolution of $(\ref{3.2})$. Therefore  $(\ref{3.2})$ has two subsolutions namely $\alpha$ and $u_1$. By \cite[theorem 3.4]{Khoi}, $(\ref{3.2})$ has a further  solution  say $u_2$,  satisfying 
$$
\max\{\alpha,u_1\}\leq u_2\leq \beta.
$$
\nd By  the maximality of $\overline{u}$, we infer that 
$$
\overline{u}(x_1) <  \overline{u}(x_2)  = u_1(x_1)  \leq u_2(x_1) \leq  \overline{u}(x_1),
$$
\nd  a contradiction. Thus  Claim \ref{rad symm}  holds true.
\vskip.3cm

\nd We set
$$
u(r) = \overline{u}(x)~\mbox{where}~r = |x|~\mbox{and}~  x \in B. 
$$
\nd and notice that  
$$
{u} \geq 0,~~{u} \ne 0,~~ u'(0)=u(R) = 0.
$$
\nd  Now, let $r \in (0,R)$ and pick $\epsilon > 0$ small such that  $r + \epsilon < R$. Remember that   $\overline{u} \in C_{0}^ {1}({\overline{B}})$ and
\begin{equation}\label{Eq u bar}
\int_B \phi(|\nabla\overline{u}|) \nabla \overline {u}\nabla v  dx =\lambda\int_B f(\overline{u}) v dx,~  v\in W_0^{1, \Phi}({B}).
\end{equation}
\nd Adapting an argument employed in \cite{santos}, consider the radially symmetric cut-off function $v_{r,\epsilon}(x) = v_{r,\epsilon}(r) $, where
$$
v_{r,\epsilon}(t) := \left\{ \begin{array}{l}
1~~ \mbox{if}~~ 0 \leq t \leq r,\\
linear~~ \mbox{if}~~ r \leq t \leq r + \epsilon,\\
0~~ \mbox{if}~~ r+\epsilon \leq t \leq R.
\end{array} \right.
$$
\nd and notice that  $v_{r,\epsilon} \in W_0^{1, \Phi}({B}) \cap Lip( {\overline{B}}    )$. Setting  $v = v_{r,\epsilon}$  in (\ref{Eq u bar})  and using the radial symmetry we get to
$$
\frac{-1}{\epsilon} \int_{r}^{r + \epsilon} t^{N-1}    
\phi(| {{u}}^{\prime}|)
 {{u}}^{\prime}~dt =
\int_0^r t^{N-1}~ \lambda
f({u})~dt + \int_{r}^{r+\epsilon} t^{N-1} \lambda f({u})\upsilon dt. \\
~~
$$
\nd Making $\epsilon \to 0$ gives
\begin{equation}\label{Eq u bar limit}
- r^{N-1} \phi(| {{u}}^{\prime}(r)|) {{u}}^{\prime}(r) =
\displaystyle \int_0^r ~ \lambda
f({u}) t^{N-1}~dt,~0 < r < R.
\end{equation}
\nd Set  
$$
\|u\|_\infty = \max \{ u(r)~|~r \in [0,R] \},
$$
\nd and choose numbers $r_0, r_1 \in [0,R)$ with $r_1 \in (r_0, R)$ such that 
$$
u(r_0) = \|u \|_{\infty}~\mbox{and}~ u(r_1) = a_1.
$$
\nd Note that
$$
u(r_0)  >  u(r_1)~\mbox{and}~ 0 \leq r_0 < r_1 < R.
$$
\begin{claim}\label{b1}
~~~~  $\| u \|_{\infty} > b_1$.
\end{claim}
\nd Indeed, assume  on the contrary that,  $u(r_0) \leq b_1$. Take $\delta > 0$ small such that 
$$
a_1 < u(r) \leq u(r_0),~~ r_0  \leq  r \leq   r_0 + \delta.
$$
\nd We have by (\ref{Eq u bar limit})
\begin{equation}\label {r0}
-r_{0}^{N-1} \phi(| {{u}}^{\prime}(r_0)|) {{u}}^{\prime}(r_0) = \displaystyle \int_{0}^{r_0} ~ \lambda
f({u}) t^{N-1}~dt,
\end{equation}
\begin{equation}\label{rr}
-r^{N-1} \phi(| {{u}}^{\prime}(r)|) {{u}}^{\prime}(r) = \displaystyle\int_{0}^{r} ~ \lambda
f({u}) t^{N-1}~dt.
\end{equation}
\nd Subtracting (\ref{rr}) minus (\ref{r0})   term by term and recalling that  $u^{\prime}(r_0) = 0$, 
$$
-r^{N-1} \phi(| {{u}}^{\prime}(r)|) {{u}}^{\prime}(r) = \int_{r_0}^{{r}} ~ \lambda
f({u}) t^{N-1}~dt,~  r_0  \leq  r \leq   r_0 + \delta.
$$
\nd Since $f \leq 0~\mbox{on}~[a_1,b_1]$, 
$$
r^{N-1} \phi(| {{u}}^{\prime}(r)|) {{u}}^{\prime}(r) \geq 0,~ r_0 \leq r \leq r_0 + \delta. 
$$ 
\nd It follows that ${{u}}^{\prime}(r) \geq 0~\mbox{for}~ r_0 \leq r \leq r_0 + \delta$. But, since $r_0$ is a global maximum on $[0,R]$, it follows that $u^{\prime} = 0~\mbox{on}~ [ r_0,  r_0 + \delta]$. By a continuation argument we get $u^{\prime} = 0~\mbox{on}~ [r_0, r_1)$ so that $u = \| u \|_{\infty}~ \mbox{on}~ [r_0, r_1]$, contradicting $u(r_0) > a_1$. As a consequence,  $\| u \|_{\infty}  >  b_1$, proving  Claim \ref{b1}.

\begin{claim}
~~ $u \in C^{2}({\cal{O}})$~where ${\cal{O}} := \{r \in (0,R)~|~ u'(r) \neq 0 \}$.
\end{claim}

\nd Of course ${\cal{O}}$  is an open set. Motivated by the left hand side of (\ref{Eq u bar limit}) consider
$$
G(z)  =   \phi(z) z,~~ z \in \r ,
$$
\nd where $z$ is set to play the role of $u^{\prime}$.  Recall that  
$$
G~\mbox{is odd},~~  G^{\prime}(z) = (\phi(z) z )^{\prime} > 0~\mbox{for}~z > 0 
$$
\nd and
$$
G(z) =  \phi(| {{u}}^{\prime}(r)|) {{u}}^{\prime}(r) =  -\frac{1}{r^{N-1}}\int_{0}^{{r}} ~ \lambda
f({u}) t^{N-1}~dt.
$$
\nd Since $\phi(z)z \in C^{1}$ and $(\phi(z)z)'\neq 0$ for $z\neq 0$, we get by applying the Inverse Function Theorem  in ${\cal{O}}$ that $z = z(r,u)$ is a $C^{1}$-function of $r$. Since $z = u^{\prime}$, the claim is proved.

\begin{claim}\label{int positive}
    $\int_{ a_{1}}^{\|u \|_{\infty}} f(s) ds > 0$.
\end{claim}

\nd Differentiating in (\ref{Eq u bar limit}) and multiplying by $u^{\prime}$ we get
$$
\big ( t^{N-1} \phi(| {{u}}^{\prime}(t)|) {{u}}^{\prime}(t)  \big)^{\prime} u^{\prime}(t) = -\lambda
f({u(t)}) u^{\prime}(t) t^{N-1},
$$
\nd and hence
$$
\big [(N-1)  t^{N-2} \phi(| u^{\prime}(t)| )u^{\prime}(t) + t^{N-1} ( \phi(| {u}^{\prime}(t)| u^{\prime}(t))^{\prime}   \big] u^{\prime}(t) = -\lambda
f({u(t)}) u^{\prime}(t) t^{N-1},
$$
\nd which gives
$$
\frac{(N-1)}{t}   \phi(| u^{\prime}|) (u^{\prime})^2 + ( \phi(| {u}^{\prime}|) u^{\prime})^{\prime}  u^{\prime} = -\lambda
f({u}) u^{\prime}.
$$
\nd Integrating from $r_0$ to $r_1$ we have
\begin{equation}\label{before}
-\Big[\int_{r_0}^{r_1}  \frac{(N-1)}{t}   \phi(| u^{\prime}|) (u^{\prime})^2  dt +\int_{r_0}^{r_1} \big[\phi(| {u}^{\prime}| )u^{\prime} \big]^{\prime}  u^{\prime}  dt \Big]=   \int _{r_0}^{r_1} \lambda f({u}) u^{\prime} dt.
\end{equation}
\nd Computing the second integral in the left hand side of (\ref{before}) we get
$$
\displaystyle \int_{r_0}^{r_1} \big [ \phi(| {{u}}^{\prime}(t)|) {{u}}^{\prime}(t)  \big]^{\prime} u^{\prime}(t)  dt  = \int_{r_0}^{r_1} [\phi(u') u^{\prime }+\phi'(u')  (u^{\prime})^{2} ]u^{''} dt~\mbox{if}~u^{\prime} > 0, 
$$
\nd and
$$
\displaystyle \int_{r_0}^{r_1} \big [ \phi(| {{u}}^{\prime}(t)|) {{u}}^{\prime}(t)  \big]^{\prime} u^{\prime}(t)  dt  = \int_{r_0}^{r_1} [\phi(-u') u^{\prime }-\phi'(-u')  (u^{\prime})^{2}]u^{''}dt~\mbox{if}~u^{\prime} < 0.
$$
\nd Making the change of variables  $s = u^{\prime}(t)$     above  we get
$$
\displaystyle \int_{r_0}^{r_1} \big [ \phi(| {{u}}^{\prime}(t)|) {{u}}^{\prime}(t)  \big]^{\prime} u^{\prime}(t)  dt  =\int_{0}^{u'(r_1)} [s\phi(|s|)]'sds.
$$
\nd Taking into (\ref{before}) we get
$$
\displaystyle \int _{\|u\|_{\infty}}^{a_1} \lambda
f({s}) ds  = 
\displaystyle -\int_{r_0}^{r_1}  \frac{(N-1)}{t}   \phi(| u^{\prime}|) (u^{\prime})^2 dt - \int_{0}^{u'(r_1)} [s\phi(|s|)]'sds < 0. 
$$
\nd On the other hand $u^{\prime}$ is not identically zero on $[r_0,r_1]$ because otherwise we would have $u(t) = u(r_0) = u(r_1)$ contradicting $a_1  < \|u\|_\infty \leq  a_2$. Therefore
$$
\int ^{\|u\|_{\infty}}_{a_1} 
f({s}) ds  > 0, 
$$
\nd proving {\bf Claim  \ref{int positive}}.
\vskip.2cm

\nd Since $f \geq 0$ on $(b_1,a_2)$ and $\|u\|_{\infty} > b_1$ it follows that 
$$
\int ^{a_2}_{a_1} 
f({s}) ds  > 0, 
$$
\nd ending the proof of  lemma \ref{thm 3.1 Kaye}.  \hfill \fbox \hsf
\vskip.2cm

\nd {\bf Case 2}~~$f(0) = 0$.
\vskip.2cm

\nd Let  $u$ be a  solution of $(\ref{1.1})$ with $a_{1} < \|u\|_\infty \leq a_{2}$. Consider a continuous function   $\widetilde{f}$ such that $\widetilde{f}(0)>0$, $\widetilde{f}(s)\geq f(s)$ if $0 \leq s \leq a_1$ and $\widetilde{f}(s)=f(s)$ if  $a_1 \leq s < \infty$. It follows that  $u$ is a subsolution of

\begin{equation}\label{3.4}
 \left\{ \begin{array}{cl}
 -\Delta_\Phi u=\lambda \widetilde{f}(u)~~\mbox{ in} ~ B \\
  \\
u\geq 0,\ \     u\in W_0^{1,\Phi}(B).
       \end{array} \right.
\end{equation}
\nd As in {\bf Case 1}  we  use $\beta(x)=a_{2}$ as a supersolution of $(\ref{3.4})$. 
\vskip.2cm

\nd Hence $(\ref{3.4})$ has a solution $\widetilde{u}$ satisfying $u\leq\widetilde{u}\leq a_2$. We now proceed as in the first part of the proof with $\widetilde{f}$ in place of $f$ to obtain 
$$
\int_{a_1}^{a_{2}}f(s)ds =\int_{a_1}^{a_{2}}\widetilde{f}(s) ds>0.
$$
\nd The proof for $a_k,~k > 2$  follows the same lines.  Theorem 1.1 is proved. \hfill \fbox \hsf
\vskip.3cm

\nd It remains to proof  lemmas \ref{lemma 3.1 Kaye} and \ref{lemma 3.2 Kaye}.  
\vskip.3cm

\nd \proof~{\bf of lemma \ref{lemma 3.1 Kaye}}. At first, using the facts that $f(0)>0$ and $(s\Phi(s))$ is strictly increasing for  $s > 0$ there is a constant $c>0$ such that
\begin{equation}\label{serrin}
\lambda f(s)+c(s\Phi(s))'\geq 0,~ s\in [0,a_k]
\end{equation}
\nd  and remember that 
\begin{equation}\label{u sol}
\int_\Omega\phi(|\nabla u|)\nabla u\nabla v dx =\lambda\int_\Omega f(u)v dx,~~ v\in W_0^{1,\Phi} .
\end{equation}
\nd Adding 
$$
c\int_\Omega (u\Phi(u))'v dx
$$
\nd to both sides of (\ref{u sol}), taking $v \geq 0$,  $0 \leq u\leq a_k$ and using (\ref{serrin}) we have
$$
 \int_\Omega\phi(|\nabla u|)\nabla u\nabla v dx +c\int_\Omega (u\Phi(u))'v dx =  \lambda\int_\Omega  \big (f(u) + c (u\Phi(u))' \big)v  dx 
   \geq 0.
$$
\nd At this point, we will  use Theorem 1.1 of \cite{Pucci}. In order  to verify condition $(1.6)$ of \cite{Pucci} set $H(s)=s\Phi'(s)-\Phi(s)$ for  $s\geq 0$ and $F(s)=cs\Phi(s)$. Note that

$$
  c\frac{s\Phi(s)}{H(s)} =  c\frac{s}{\frac{s\Phi'(s)}{\Phi(s)}-1}     
   \leq \frac{cs}{\gamma_1-1}, 
$$
\nd where in the last inequality we used $(\ref{Delta_2})$. 
\vskip.2cm

\nd By the inequality above  choose $\delta>0$ such that $c\frac{s\Phi(s)}{H(s)}\leq 1$ for $s\in (0,\delta)$. Since  $H^{-1}$ is strictly increasing, we infer that $H^{-1}(cs\Phi(s))\leq s$ for $s\in (0,\delta)$, from which condition $(1.6)$ of \cite{Pucci} follows.  
\vskip.1cm

\nd This ends the proof of lemma \ref{lemma 3.1 Kaye}.   \hfill \fbox \hsf
\vskip.2cm

\nd \proof~ {\bf of lemma \ref{lemma 3.2 Kaye}}~  Of course $\beta$ is an upper-solution of  (\ref{3.2}). To deal with $\alpha$ define 
$$
v_n(x)=n\min\{u(x),\frac{1}{n}\}~\mbox{ for}~ x\in\Omega,~n \geq 1~\mbox{is an integer}.
$$
\nd   Notice that $\nabla u \cdot \nabla v_n\geq 0$ and by the very definition, $v_n$ converges to $1$, pointwisely  in $\Omega$. Take $w\geq 0$, $w\in C_0^\infty(B)$ and note that $wv_n\in W_0^{1,\Phi}(\Omega)$. This implies 
$$
\int_\Omega \phi(|\nabla u|)\nabla u\nabla (wv_n) dx =\lambda\int_\Omega f(u)wv_n dx
$$
\nd which gives 

$$\int_\Omega w\phi(|\nabla u|)\nabla u\nabla v_n dx +\int_\Omega v_n\phi(|\nabla u|)\nabla u\nabla w dx=\lambda\int_\Omega f(u)wv_n dx$$

\nd We observe that 
$$
0\leq w\phi(|\nabla u|)\nabla u\nabla v_n\leq w\phi(|\nabla u|)|\nabla u|^2~\mbox{and}~  0\leq v_n\leq 1.
$$

\nd By the Lebesgue Theorem we infer that
\begin{eqnarray}
  \int_B \phi(|\nabla\alpha|)\nabla \alpha\nabla w dx&=&  \int_\Omega \phi(|\nabla u|)\nabla u\nabla w  dx    \nonumber \\
   &=& \lim\int_\Omega v_n\phi(|\nabla u|)\nabla u\nabla w dx\nonumber \\
   &=& \lim\int_\Omega (\lambda f(u)wv_n-w\phi(|\nabla u|)\nabla u\nabla v_n) dx \nonumber \\
   &\leq& \lambda\int_\Omega f(u)w dx\nonumber \\
   &\leq& \lambda\int_B f(\alpha)w dx. \nonumber
\end{eqnarray}
\nd  This ends the proof of lemma \ref{lemma 3.2 Kaye}.  \hfill \fbox \hsf

\section {Appendix }

\nd  In this section we present for the sake of completeness  several rather simple results and remarks which were employed in the body of the paper.

\begin{prop}\label{equivalencia-phi34}

  If $\phi$ sartisfies  $(\phi_3)-(\phi_4)$ then $(\ref{relacao a_1- a_2})$ holds true. 
\end{prop}

\nd \proof~ Indeed by  $(\phi_3)$,
$$
    \displaystyle \sum_{i,j=1}^N\frac{\partial \alpha_j}{\partial\eta_i}(\eta) \delta_{i,j} \geq\Gamma_1~ \phi(|\eta|)|\xi|^2,
$$
\nd where $\delta_{i,j}$ is the Kronecker symbol. Thus
$$
    \displaystyle \sum_{i,j=1}^N   \frac{\partial  [{\phi(|\eta|) \eta_{j} }]} {\partial\eta_i} \delta_{i,j} \geq\Gamma_1~ \phi(|\eta|)|\xi|^2.
$$
\nd Take $t > 0$, $\eta=(t,0,...,0)$  and  $\xi=(1,0,...,0)$. Then
$$
\begin{array}{rcc}
   \displaystyle \sum_{i,j=1}^N   \frac{\partial  [{\phi(|\eta|) \eta_{j} }]} {\partial\eta_i} \delta_{i,j} =
\frac{d (t \phi(t) )}{dt}\\

\geq \Gamma_1~ \phi(t).
\end{array}
$$
\nd Therefore 
$$
\frac{ ( t \phi(t) )^{\prime} }{\phi(t)} \geq \Gamma_1,~t > 0.
$$
\nd On the other hand, assuming $(\phi_4)$ we have
$$
   \displaystyle \sum_{i,j=1}^N   \frac{\partial  [{\phi(|\eta|) \eta_{j} }]} {\partial\eta_i} \delta_{i,j} 
\leq\Gamma_2~ \phi(|\eta|),
$$
\nd Take $t > 0$, $\eta=(t,0,...,0)$  and  $\xi=(1,0,...,0)$. Arguing as above we find 
$$
\frac{(t\phi(t))^{\prime}}{   \phi(t)    }\leq \Gamma_2,~t > 0.
$$
\nd This ends the proof of proposition \ref{equivalencia-phi34}.   \hfill \fbox \hsf
\vskip.1cm

\begin{rmk}\label{Aa} {\rm Verification of  {\bf (iii)} in remark (\ref{obss})}.
\vskip.1cm
\nd By  $(\ref{relacao a_1- a_2})$  we have
$$
\Gamma_{1} \phi(s) \leq (s\phi(s))^{\prime} \leq \Gamma_{2} \phi(s),~~s > 0.
$$
\nd Multiplying by $s$ and integrating  from $0$ to $t$ we have
$$
\Gamma_{1} \Phi(t) \leq  t^2 \phi(t) - \Phi(t) \leq \Gamma_{2} \Phi(t),~~t > 0.
$$
\nd As a consequence,  
$$
(\Gamma_{1} + 1) \Phi(t) \leq  t\Phi^{\prime}(t) \leq (\Gamma_{2} + 1) \Phi(t),~~t > 0,
$$
\nd showing $(\ref{Delta_2})$.
\end{rmk}

\begin{rmk} {\bf (On example \ref{ex 1}) }.
\vskip.2cm

\nd Let $\phi(t)=2\gamma (1+t^2)^{\gamma-1}$ with $\gamma>\frac{1}{2}$.  Then $\Phi(t)=(1+t^2)^\gamma - 1$.
\vskip.1cm

\nd Differentiating in the expression of $\phi$ we get
$$
\phi'(t)=4\gamma(\gamma-1) (1+t^2)^{\gamma-2}t.
$$
\nd It follows that
$$
\frac{(t\phi(t))'}{\phi(t)}=1+2(\gamma-1)\frac{t^2}{1+t^2}.
$$
\nd and so
$$
\min\{1,2\gamma-1\}\leq\frac{(t\phi(t))'}{\phi(t)}\leq\max\{1,2\gamma-1\}.
$$
\nd By  proposition \ref{equivalencia-phi34}, $\phi$ satisfies $(\phi_3)-(\phi_4)$.  It follows that $\phi$ satisfies   $(\phi_i),~i = 1, \cdots,4$.   \hfill \fbox \hsf
\end{rmk}

\begin{rmk} {\bf (On example \ref{ex 2}) }.
\vskip.2cm

 \nd { \rm Consider
$$
\phi(t) =\frac{pt^{p-2}(1+t)\ln(1+t)+t^{p-1}}{1+t},~t > 0.
$$
\nd Then  
$$
\Phi(t)=t^p\ln(1+t).
$$
\nd By computing, we get
$$
(t\phi(t))'=t^{p-2}\left[p(p-1)\ln(1+t)+\frac{2p t}{1+t}-\frac{t^2}{(1+t)^2}\right]
$$
\nd so that 
$$
\frac{(t\phi(t))'}{\phi(t)}=\frac{2 p (1 + t)t-t^2 + p(p-1) (1 + t)^2 \ln(1 + t)}
{(1 + t) (t + p (1 + t) \ln(1 + t))},
$$
\nd which is a decreasing function. Moreover, 
$$
\lim_{t\rightarrow\infty}\frac{(t\phi(t))'}{\phi(t)}=p-1~~\mbox{and}~~\lim_{t\rightarrow 0}\frac{(t\phi(t))'}{\phi(t)}=p
$$
\nd and so
$$
p-1\leq\frac{(t\phi(t))'}{\phi(t)}\leq p, 
$$
\nd By proposition  \ref{equivalencia-phi34}, $\phi$ satisfies $(\phi_3)-(\phi_4)$.\hfill \fbox \hsf }
\end{rmk}

\nd We refer the reader to  \cite{Fukagai} and references therein for the  lemma below whose proof is elementary.
\begin{lem}\label{lema_naru}
       Assume that  $\phi$ satisfies  $(\phi_1)-(\phi_3)$.
        Set
         $$
         \zeta_0(t)=\min\{t^{\gamma_1},t^{\gamma_2}\},~~~ \zeta_1(t)=\max\{t^{\gamma_1},t^{\gamma_2}\},~~ t\geq 0.
        $$
  \nd Then  $\Phi$ satisfies
       $$
            \zeta_0(t)\Phi(\rho)\leq\Phi(\rho t)\leq \zeta_1(t)\Phi(\rho),~~ \rho, t> 0,
        $$
$$
\zeta_0(\|u\|_{\Phi})\leq\int_\Omega\Phi(u)dx\leq \zeta_1(\|u\|_{\Phi}),~ u\in L_{\Phi}(\Omega).
 $$
\end{lem}

\begin{flushright}
{\small M. L. Carvalho}\\ 
{\small J. V. Goncalves}\\
{\small K. O. Silva}\\
\smallskip
  \scriptsize{Universidade Federal de Goi\'as\\
   Instituto de Matem\'atica e Estat\'istica\\
   74001-970 Goi\^ania, GO - Brasil}
\end{flushright}

\end{document}